\documentclass[fleqn]{llncs}
\usepackage{pkfikf}
\usepackage{version}

\title{Partial Komori Fields and\\ Imperative Komori Fields}
\author{J.A. Bergstra \and C.A. Middelburg}
\institute{Informatics Institute, Faculty of Science,
           University of Amsterdam, \\
           Science Park~107, 1098~XG Amsterdam, the Netherlands \\
           \email{J.A.Bergstra@uva.nl,C.A.Middelburg@uva.nl}}

\begin{document}

\maketitle

\begin{abstract}
% 94 + 4 %
This paper is concerned with the status of $1 \mdiv 0$ and ways to deal
with it.
These matters are treated in the setting of Komori fields, also known as
non-trivial cancellation meadows.
Different viewpoints on the status of $1 \mdiv 0$ exist in mathematics
and theoretical computer science.
We give a simple account of how mathematicians deal with $1 \mdiv 0$ in
which a customary convention among mathematicians plays a prominent
part, and we make plausible that a convincing account, starting from the
popular computer science viewpoint that $1 \mdiv 0$ is undefined, by
means of some logic of partial functions is not attainable.
\begin{classcode}
12E99, 12L99.
\end{classcode}
\end{abstract}

\section{Introduction}
\label{sect-introduction}

In~\cite{BT07a}, meadows are proposed as alternatives for fields with a
purely equational specification.
A meadow is a commutative ring with identity and a total multiplicative
inverse operation satisfying two equations which imply that the
multiplicative inverse of zero is zero.
Thus, meadows are total algebras.
Recently, we found in~\cite{Ono83a} that meadows were already introduced
by Komori~\cite{Kom75a} in a report from 1975, where they go by the name
of \emph{desirable pseudo-fields}.
This finding induced us to propose the name \emph{Komori field} for a
meadow satisfying $0 \neq 1$ and $x \neq 0 \Implies x \mmul x\minv = 1$
(see~\cite{BM09h}).
The prime example of Komori fields is the field of rational numbers with
the multiplicative inverse operation made total by imposing that the
multiplicative inverse of zero is zero.

In~\cite{BM09g}, we renamed meadows to inversive meadows and introduced
divisive meadows as inversive meadows with the multiplicative inverse
operation replaced by a division operation.
We also introduced simple constructions of variants of inversive and
divisive meadows with a partial multiplicative or division operation.
Moreover, we took a survey of logics of partial functions, but did not
succeed in determining their adequacy for reasoning about the partial
variants of meadows.
This made us reflect on the way in which mathematicians deal with
$1 \mdiv 0$ and how that way relates to the viewpoint, reflected in the
partial variants of meadows, that $1 \mdiv 0$ is undefined.

In the current paper, we give a simple account of how mathematicians
deal with $1 \mdiv 0$ in mathematical works.
Dominating in this account is the concept of an imperative Komori field,
a concept in which a customary convention among mathematicians plays a
prominent part.
We also make plausible that a convincing account, starting from the
usual viewpoint of theoretical computer scientists that $1 \mdiv 0$ is
undefined, by means of some logic of partial functions is not
\linebreak[2] attainable.

It is quite usual that neither the division operator nor the
multiplicative inverse operator is included in the signature of number
systems such as the field of rational numbers and the field of real
numbers.
However, the abundant use of the division operator in mathematical
practice makes it very reasonable to include the division operator, or
alternatively the multiplicative inverse operator, in the signature.
It appears that excluding both of them creates more difficulties than
that it solves.
At the least, the problem of division by zero cannot be avoided by
excluding $1 \mdiv 0$ from being written.

This paper is organized as follows.
First, we discuss the main prevailing viewpoints on the status of $1
\mdiv 0$ in mathematics and theoretical computer science
(Section~\ref{sect-viewpoints-div-by-zero}).
Next, we give a brief summary of meadows and Komori fields
(Section~\ref{sect-KF}).
After that, we introduce partial Komori fields and imperative Komori
fields (Section~\ref{sect-PKF-IKF}) and discuss the convention that is
involved in imperative Komori fields (Section~\ref{sect-RDC-discussed}).
Then, we make plausible the inadequacy of logics of partial functions
for a convincing account of how mathematicians deal with $1 \mdiv 0$
(Section~\ref{sect-inadequacy-LPF}).
Finally, we make some concluding remarks
(Section~\ref{sect-conclusions}).

\section{Viewpoints on the Status of $1 \mdiv 0$}
\label{sect-viewpoints-div-by-zero}

In this section, we shortly discuss two prevailing viewpoints on the
status of $1 \mdiv 0$ in mathematics and one prevailing viewpoint on
the status of $1 \mdiv 0$ in theoretical computer science.
To our knowledge, the viewpoints in question are the main prevailing
viewpoints.
We take the case of the rational numbers, the case of the real numbers
being essentially the same.

One prevailing viewpoint in mathematics is that $1 \mdiv 0$ has no
meaning because $1$ cannot be divided by $0$.
The argumentation for this viewpoint rests on the fact that there is no
rational number $z$ such that $0 \mmul z = 1$.
Moreover, in mathematics, syntax is not prior to semantics and posing
the question ``what is $1 \mdiv 0$'' is not justified by the mere
existence of $1 \mdiv 0$ as a syntactic object.
Given the fact that there is no rational number that mathematicians
intend to denote by $1 \mdiv 0$, this means that there is no need to
assign a meaning to $1 \mdiv 0$.

Another prevailing viewpoint in mathematics is that the use of
$1 \mdiv 0$ is simply disallowed because the intention to divide $1$ by
$0$ is non-existent in mathematical practice.
This viewpoint can be regarded as a liberal form of the previous one:
the rejection of the possibility that $1 \mdiv 0$ has a meaning is
circumvented by disallowing the use of $1 \mdiv 0$.
Admitting that $1 \mdiv 0$ has a meaning, such as $0$ or ``undefined'',
is consistent with this viewpoint.

The prevailing viewpoint in theoretical computer science is that the
meaning of $1 \mdiv 0$ is ``undefined'' because division is a partial
function.
Division is identified as a partial function because there is no
rational number $z$ such that $0 \mmul z = 1$.
This viewpoint presupposes that the use of $1 \mdiv 0$ should be
allowed, for otherwise assigning a meaning to $1 \mdiv 0$ does not make
sense.
Although this viewpoint is more liberal than the previous one, it is
remote from ordinary mathematical practice.
This will be illustrated in Section~\ref{sect-inadequacy-LPF}.

The first of the two prevailing viewpoints in mathematics discussed
above only leaves room for very informal concepts of expression,
calculation, proof, substitution, etc.
For that reason, we refrain from considering that viewpoint any further.

The prevailing viewpoint in theoretical computer science corresponds to
two of the partial meadows of rational numbers obtained from the
inversive and divisive meadows of rational numbers by a simple
construction in~\cite{BM09g}.
Because $0 \neq 1$ and $x \neq 0 \Implies x \mmul x\minv = 1$ are
satisfied, the latter are also called Komori fields of rational numbers
and the former are also called partial Komori fields of rational
numbers.

The prevailing viewpoint in mathematics considered further in this paper
corresponds to the inversive and divisive meadows of rational numbers
together with an imperative about the use of the multiplicative inverse
operator and division operator, respectively.
These combinations are called imperative Komori fields.

\section{Meadows and Komori Fields}
\label{sect-KF}

In this section, we give a brief summary of meadows and Komori fields.

The signature of inversive meadows consists of the following constants
and operators:
\begin{itemize}
\item
the constants $0$ and $1$;
\item
the binary \emph{addition} operator ${} +$ {};
\item
the binary \emph{multiplication} operator ${} \mmul {}$;
\item
the unary \emph{additive inverse} operator $- {}$;
\item
the unary \emph{multiplicative inverse} operator ${}\minv$.
\end{itemize}
We use infix notation for the binary operators, prefix notation for the
unary operator $- {}$, and postfix notation for the unary operator
${}\minv$.
Moreover, we use the usual precedence convention to reduce the need for
parentheses.

The set of all terms over the signature of inversive meadows constitutes
the \emph{inversive meadow notation}.

An inversive meadow is an algebra over the signature of inversive
meadows that satisfies the equations given in
Tables~\ref{eqns-commutative-ring} and~\ref{eqns-add-inversive-meadow}.%
\begin{table}[!t]
\caption{Axioms of a commutative ring with identity}
\label{eqns-commutative-ring}
\begin{eqntbl}
\begin{eqncol}
(x + y) + z = x + (y + z)                                             \\
x + y = y + x                                                         \\
x + 0 = x                                                             \\
x + (-x) = 0
\end{eqncol}
\qquad\quad
\begin{eqncol}
(x \mmul y) \mmul z = x \mmul (y \mmul z)                             \\
x \mmul y = y \mmul x                                                 \\
x \mmul 1 = x                                                         \\
x \mmul (y + z) = x \mmul y + x \mmul z
\end{eqncol}
\end{eqntbl}
\end{table}
\begin{table}[!t]
\caption{Additional axioms for an inversive meadow}
\label{eqns-add-inversive-meadow}
\begin{eqntbl}
\begin{eqncol}
(x\minv)\minv = x                                                   \\
x \mmul (x \mmul x\minv) = x
\end{eqncol}
\end{eqntbl}
\end{table}
The equations given in Table~\ref{eqns-commutative-ring} are the axioms
of a commutative ring with identity.
From the equations given in Tables~\ref{eqns-commutative-ring}
and~\ref{eqns-add-inversive-meadow}, the equation $0\minv = 0$ is
derivable.

Henceforth, we will write $\sigimd$ for the signature of inversive
meadows and $\eqnsimd$ for the set of axioms for inversive meadows.

A \emph{inversive Komori field} is an inversive meadow that satisfies
the \emph{separation axiom} $0 \neq 1$ and
the \emph{general inverse law} $x \neq 0 \Implies x \mmul x\minv = 1$.

The inversive Komori field that we are most interested in is $\Ratzi$,
the inversive Komori field of rational numbers:%
\footnote
{We write $I(\Sigma,E)$ for the initial algebra among the algebras
 over the signature $\Sigma$ that satisfy the equations $E$
 (see e.g.~\cite{BT87a}).}
\pagebreak[2]
\begin{ldispl}
\Ratzi =
I(\sigimd,\eqnsimd \union
  \set{(1 + x^2 + y^2) \mmul (1 + x^2 + y^2)\minv = 1})\;.
\end{ldispl}
$\Ratzi$ differs from the field of rational numbers only in that the
multiplicative inverse of zero is zero.

The signature of divisive meadows is the signature of inversive meadows
with the unary multiplicative inverse operator ${}\minv$ replaced by:
\begin{itemize}
\item
the binary \emph{division} operator ${} \mdiv {}$.
\end{itemize}

The set of all terms over the signature of divisive meadows constitutes
the \emph{divisive meadow notation}.

A divisive meadow is an algebra over the signature of divisive
meadows that satisfies the equations given in
Tables~\ref{eqns-commutative-ring} and~\ref{eqns-add-divisive-meadow}.%
\begin{table}[!t]
\caption{Additional axioms for a divisive meadow}
\label{eqns-add-divisive-meadow}
\begin{eqntbl}
\begin{eqncol}
1 \mdiv (1 \mdiv x) = x                                               \\
(x \mmul x) \mdiv x = x                                               \\
x \mdiv y = x \mmul (1 \mdiv y)
\end{eqncol}
\end{eqntbl}
\end{table}

Henceforth, we will write $\sigdmd$ for the signature of divisive
meadows and $\eqnsdmd$ for the set of axioms for divisive meadows.

A \emph{divisive Komori field} is an divisive meadow that satisfies
the \emph{separation axiom} $0 \neq 1$ and
the \emph{general division law} $x \neq 0 \Implies x \mdiv x = 1$.

The divisive Komori field that we are most interested in is $\Ratzd$,
the divisive Komori field of rational numbers:
\begin{ldispl}
\Ratzd =
I(\sigdmd,\eqnsdmd \union
  \set{(1 + x^2 + y^2) \mdiv (1 + x^2 + y^2) = 1})\;.
\end{ldispl}
$\Ratzd$ differs from $\Ratzi$ only in that the multiplicative inverse
operation is replaced by a division operation such that
$x \mdiv y = x \mmul y\minv$.

\section{Partial Komori Fields and Imperative Komori Fields}
\label{sect-PKF-IKF}

In this section, we introduce partial inversive and divisive Komori
fields and imperative inversive and divisive Komori fields.

Let $\cK_i$ be an inversive Komori field.
Then it make sense to construct one partial inversive Komori field from
$\cK_i$:
\begin{itemize}
\item
$0\minv \punch \cK_i$ is the partial algebra that is obtained from
$\cK_i$ by making $0\minv$ undefined.
\end{itemize}
Let $\cK_d$ be a divisive Komori field.
Then it make sense to construct two partial divisive Komori fields from
$\cK_d$:
\begin{itemize}
\item
$\Quant \mdiv 0 \punch \cK_d$ is the partial algebra that is obtained
from $\cK_d$ by making\linebreak[2] $q \mdiv 0$ undefined for all $q$ in
the domain of $\cK_d$;
\item
$(\Quant \diff \set{0}) \mdiv 0 \punch \cK_d$ is the partial algebra
that is obtained from $\cK_d$ by\linebreak[2] making $q \mdiv 0$
undefined for all $q$ in the domain of $\cK_d$ different from $0$.
\end{itemize}
$(\Quant \diff \set{0}) \mdiv 0 \punch \cK_d$ expresses a view on the
partiality of division by zero that cannot be expressed if only
multiplicative inverse is available.

The partial Komori field constructions are special cases of a
more general partial algebra construction for which we have coined the
term \emph{punching} in~\cite{BM09g}.
Presenting the details of the general construction is outside the scope
of the current paper.

The partial Komori fields that we are most interested in are the ones
that can be obtained from $\Ratzi$ and $\Ratzd$ by means of the partial
Komori field constructions introduced above.
It yields three partial Komori fields of rational numbers:
\begin{ldispl}
\begin{eqncol}
0\minv \punch \Ratzi\;,
\qquad
\Quant \mdiv 0 \punch \Ratzd\;,
\qquad
(\Quant \diff \set{0}) \mdiv 0 \punch \Ratzd\;.
\end{eqncol}
\end{ldispl}

Notice that these partial algebras have been obtained by means of the
well-known initial algebra construction and a straightforward partial
algebra construction.
This implies that only equational logic for total algebras has been used
as a logical tool for their construction.

The first two of the partial Komori fields of rational numbers
introduced above correspond most closely to the prevailing viewpoint on
the status of $1 \mdiv 0$ in theoretical computer science that is
mentioned in Section~\ref{sect-viewpoints-div-by-zero}.
In the sequel, we will focus on $\Quant \mdiv 0 \punch \Ratzd$ because
the divisive notation is used more often than the inversive notation.

An imperative Komori field of rational numbers is a Komori field of
rational numbers together with an imperative to comply with a very
strong convention with regard to the use of the multiplicative inverse
or division operator.

Like with the partial Komori field of rational numbers, we introduce
three imperative Komori fields of rational numbers:
\begin{itemize}
\item
$0\minv \imper \Ratzi$ is $\Ratzi$ together with the imperative to
comply with the convention that $q\minv$ is not used with $q = 0$;
\item
$\Quant \mdiv 0 \imper \Ratzd$ is $\Ratzd$ together with the
imperative to comply with the convention that $p \mdiv q$ is not used
with $q = 0$;
\item
$(\Quant \diff \set{0}) \mdiv 0 \imper \Ratzd$ is $\Ratzd$ together
with the imperative to comply with the convention that $p \mdiv q$ is
not used with $q = 0$ if $p \neq 0$.
\end{itemize}
The conventions are called the \emph{relevant inversive convention},
the \emph{relevant division convention} and
the \emph{liberal relevant division convention}, respectively.

The conventions are very strong in the settings in which they must be
complied with.
For example, the relevant division convention is not complied with if
the question ``what is $1 \mdiv 0$'' is posed.
Using $1 \mdiv 0$ is disallowed, although we know that $1 \mdiv 0 = 0$
in~$\Ratzd$.

The first two of the imperative Komori fields of rational numbers
introduced above correspond most closely to the second of the two
prevailing viewpoints on the status of $1 \mdiv 0$ in mathematics
that are mentioned in Section~\ref{sect-viewpoints-div-by-zero}.
In the sequel, we will focus on $\Quant \mdiv 0 \imper \Ratzd$ because
the divisive notation is used more often than the inversive notation.

Komori fields go by the name of non-trivial cancellation meadows in
previous work (see e.g.~\cite{BM09g}).
Therefore, we take the names partial non-trivial cancellation meadow
and imperative non-trivial cancellation meadow as alternatives for
partial Komori field and imperative Komori field, respectively.

\section{Discussion on the Relevant Division Convention}
\label{sect-RDC-discussed}

In this section, we discuss the relevant division convention, i.e.\ the
convention that plays a prominent part in imperative Komori fields.

The existence of the relevant division convention can be explained by
assuming a context in which two phases are distinguished: a definition
phase and a working phase.
A mathematician experiences these phases in this order.
In the definition phase, the status of $1 \mdiv 0$ is dealt with
thoroughly so as to do away with the necessity of reflection upon it
later on.
As a result, $\Ratzd$ and the relevant division convention come up.
In the working phase, $\Ratzd$ is simply used in compliance with the
relevant division convention when producing mathematical texts.
Questions relating to $1 \mdiv 0$ are understood as being part of the
definition phase, and thus taken out of mathematical practice.
This corresponds to a large extent with how mathematicians work.

In the two phase context outlined above, the definition phase can be
made formal and logical whereas the results of this can be kept out of
the working phase.
Indeed, in mathematical practice, we find a world where logic does not
apply and where validity of work is not determined by the intricate
details of a very specific formal definition but rather by the consensus
obtained by a group of readers and writers.

Whether a mathematical text, including definitions, questions, answers,
conjectures and proofs, complies with the relevant division convention
is a judgement that depends on the mathematical knowledge of the reader
and writer.
For example, $\Forall{x}{(x^2 + 1) \mdiv (x^2 + 1) = 1}$ complies with
the relevant division convention because the reader and writer of it
both know that $\Forall{x}{x^2 + 1 \neq 0}$.

Whether a mathematical text complies with the relevant division
convention may be judged differently even with sufficient mathematical
knowledge.
This is illustrated by the following mathematical text, where $>$ is the
usual ordering on the set of rational numbers:
\begin{quote}
\textbf{Theorem.}\,\,
If $p \mdiv q = 7$ then
$\displaystyle \frac{q^2 + p \mdiv q - 7}{q^4 + 1} > 0$.

\textit{Proof.}\,\,
Because $q^4 + 1 > 0$, it is sufficient to show that
$q^2 + p \mdiv q - 7 > 0$.
It follows from $p \mdiv q = 7$ that $q^2 + p \mdiv q - 7 = q^2$,
and $q^2 > 0$ because $q \neq 0$ (as $p \mdiv q = 7$). \qed
\end{quote}
Reading from left to right, it cannot be that first $p \mdiv q$ is used
while knowing that $q \neq 0$ and that later on $q \neq 0$ is inferred
from the earlier use of $p \mdiv q$.
However, it might be said that the first occurrence of the text fragment
$p \mdiv q = 7$ introduces the knowledge
that $q \neq 0$ at the right time, i.e.\ only after it has been entirely
read.

The possibility of different judgements with sufficient mathematical
knowledge looks to be attributable to the lack of a structure theory of
mathematical text.
However, with a formal structure theory of mathematical text, we still
have to deal with the fact that compliance with the relevant division
convention is undecidable.

The imperative to comply with the relevant division conventions boils
down to the disallowance of the use of $1 \mdiv 0$,
$1 \mdiv (1 + (-1))$, etcetera in mathematical text.
The usual explanation for this is the non-existence of a $z$ such that
$0 \mmul z = 1$.
This makes the legality of $1 \mdiv 0$ comparable to the legality of
$\sum_{m = 1}^\infty 1 \mdiv m$, because of the non-existence of the
limit of $\tup{\sum_{m = 1}^{n+1} 1 \mdiv m}_{n \in \Nat}$.
However, a mathematical text may contain the statement
``$\sum_{m = 1}^\infty 1 \mdiv m$ is divergent''.
That is, the use of $\sum_{m = 1}^\infty 1 \mdiv m$ is not disallowed.
So the fact that there is no rational number that mathematicians intend
to denote by an expression does not always lead to the disallowance of
its use.

In the case of $1 \mdiv 0$, there is no rational number that
mathematicians intend to denote by $1 \mdiv 0$, there is no real number
that mathematicians intend to denote by $1 \mdiv 0$, there is no complex
number that mathematicians intend to denote by $1 \mdiv 0$, etcetera.
A slightly different situation arises with $\sqrt{2}$: there is no
rational number that mathematicians intend to denote by $\sqrt{2}$, but
there is a real number that mathematicians intend to denote by
$\sqrt{2}$.
It is plausible that the relevant division convention has emerged
because there is no well-known extension of the field of rational
numbers with a number that mathematicians intend to denote by
$1 \mdiv 0$.

\section{Partial Komori Fields and Logics of Partial Functions}
\label{sect-inadequacy-LPF}

In this section, we adduce arguments in support of the statement that
partial Komori fields together with logics of partial functions do not
quite explain how mathematicians deal with $1 \mdiv 0$ in mathematical
works.
It needs no explaining that a real proof of this statement is out of the
question.
However, we do not preclude the possibility that more solid arguments
exist.
Moreover, as it stands, it is possible that our argumentation leaves
room for controversy.

In the setting of a logic of partial functions, there may be terms whose
value is undefined.
Such terms are called non-denoting terms.
Moreover, often three truth values, corresponding to true, false and
neither-true-nor-false, are considered.
These truth values are denoted by $\True$, $\False$, and $\Undef$,
respectively.

In logics of partial functions, three different kinds of equality are
found (see e.g.~\cite{MR91a,BM09g}).
They only differ in their treatment of non-denoting terms:
\begin{itemize}
\item
\emph{weak equality}: if either $t$ or $t'$ is non-denoting, then the
truth value of $t = t'$ is $\Undef$;
\item
\emph{strong equality}: if either $t$ or $t'$ is non-denoting, then the
truth value of $t = t'$ is $\True$ whenever both $t$ and $t'$ are
non-denoting and $\False$ otherwise;
\item
\emph{existential equality}: if either $t$ or $t'$ is non-denoting, then
the truth value of $t = t'$ is $\False$.
\end{itemize}

With strong equality, the truth value of $1 \mdiv 0 = 1 \mdiv 0 + 1$ is
$\True$.
This does not at all fit in with mathematical practice.
With existential equality, the truth value of $1 \mdiv 0 = 1 \mdiv 0$ is
$\False$.
This does not at all fit in with mathematical practice as well.
Weak equality is close to mathematical practice: the truth value of an
equation is neither $\True$ nor $\False$ if a term of the form
$p \mdiv q$ with $q = 0$ occurs in it.

This means that the classical logical connectives and quantifiers must
be extended to the three-valued case.
Many ways of extending them must be considered uninteresting for a logic
of partial functions because they lack an interpretation of the third
truth value that fits in with its origin: dealing with non-denoting
terms.
If those ways are excluded, only four ways to extend the classical
logical connectives to the three-valued case remain
(see e.g.~\cite{BBR95a}).
Three of them are well-known: they lead to Bochvar's strict
connectives~\cite{Boc39a}, McCarthy's sequential
connectives~\cite{McC63a}, and Kleene's monotonic
connectives~\cite{Kle38a}.
The fourth way leads to McCarthy's sequential connectives with the role
of the operands of the binary connectives reversed.

In mathematical practice, the truth value of
$\Forall{x}{x \neq 0 \Implies x \mdiv x = 1}$ is considered $\True$.
Therefore, the truth value of $0 \neq 0 \Implies 0 \mdiv 0 = 1$ is
$\True$ as well.
With Bochvar's connectives, the truth value of this formula is $\Undef$.
With McCarthy's or Kleene's connectives the truth value of this formula
is $\True$.
However, unlike with Kleene's connectives, the truth value of the
seemingly equivalent $0 \mdiv 0 = 1 \Or 0 = 0$ is $\Undef$ with
McCarthy's connectives.
Because this agrees with mathematical practice, McCarthy's connectives
are closest to mathematical practice.

The conjunction and disjunction connectives of Bochvar and the
conjunction and disjunction connectives of Kleene have natural
generalizations to quantifiers, which are called Bochvar's quantifiers
and Kleene's quantifiers, respectively.
Both Bochvar's quantifiers and Kleene's quantifiers can be considered
generalizations of the conjunction and disjunction connectives of
McCarthy.%
\footnote
{In~\cite{KTB91a}, Bochvar's quantifiers are called McCarthy's
 quantifiers, but McCarthy combines his connectives with Kleene's
 quantifiers (see e.g.~\cite{Kle38a}).}

With Kleene's quantifiers, the truth value of
$\Forall{x}{x \mdiv x = 1}$ is $\Undef$ and the truth value of
$\Exists{x}{x \mdiv x = 1}$ is $\True$.
The latter does not at all fit in with mathematical practice.
Bochvar's quantifiers are close to mathematical practice: the truth
value of a quantified formula is neither $\True$ nor $\False$ if a term
of the form $p \mdiv q$ with $q = 0$ occurs in it.

What precedes suggest that mathematical practice is best approximated by
a logic of partial functions with weak equality, McCarthy's connectives
and Bochvar's quantifiers.
We call this logic the \emph{logic of partial meadows}, abbreviated
$\LPMd$.
% and write $\Intlpmd{\phi}$, where $\phi$ is a sentence, for the
% truth value of $\phi$ according to the semantics of $\LPMd$.

In order to explain how mathematicians deal with $1 \mdiv 0$ in
mathematical works, we still need the convention that a sentence is not
used if its truth value is neither $\True$ nor $\False$.
We call this convention the \emph{two-valued logic convention}.

$\LPMd$ together with the imperative to comply with the two-valued logic
convention gets us quite far in explaining how mathematicians deal with
$1 \mdiv 0$ in mathematical works.
However, in this setting, not only the truth value of
$0 \neq 0 \Implies 0 \mdiv 0 = 1$ is $\True$, but also the truth value
of $0 = 0 \Or 0 \mdiv 0 = 1$ is $\True$.\linebreak[2]
In our view, the latter does not fit in with how mathematicians deal
with $1 \mdiv 0$ in mathematical works.
Hence, we conclude that $\LPMd$, even together with the imperative to
comply with the two-valued logic convention, fails to provide a
convincing account of how mathematicians deal with $1 \mdiv 0$ in
mathematical works.

\section{Conclusions}
\label{sect-conclusions}

We have given a simple account of how mathematicians deal with
$1 \mdiv 0$ in mathematical works.
Dominating in this account is the concept of an imperative Komori field.
The concept of an imperative Komori field is a special case of the more
general concept of an \emph{imperative algebra}, i.e.\ an algebra
together with the imperative to comply with some convention about its
use.
An example of an imperative algebra is imperative stack algebra: stack
algebra, whose signature consists of $\nm{empty}$, $\nm{push}$,
$\nm{pop}$ and $\nm{top}$, together with the imperative to comply with
the convention that $\nm{top}(s)$ is not used with $s = \nm{empty}$.

Moreover, we have argued that a logic of partial functions with weak
equality, McCarthy's connectives and Bochvar's quantifiers, together
with the imperative to comply with the convention that sentences whose
truth value is neither $\True$ nor $\False$ are not used, approximates
mathematical practice best, but after all fails to provide a convincing
account of how mathematicians deal with $1 \mdiv 0$ in math\-ematical
works.

\bibliographystyle{spmpsci}
\bibliography{MD}

% \par \vfill \par \noindent DRAFT of \today

\end{document}